\documentclass[11pt]{article}

\usepackage{amsfonts, amsbsy, amsmath, amssymb}

\renewcommand{\to}{\rightarrow}

\newcommand{\Z}{\mathbb{Z}}

\newcommand{\id}{\mbox{id}}

\newtheorem{theorem}{Theorem}[section]
\newtheorem{corollary}[theorem]{Corollary}
\newtheorem{lemma}[theorem]{Lemma}
\newtheorem{proposition}[theorem]{Proposition}

\newtheorem{definition}[theorem]{Definition}

\newtheorem{example}[theorem]{Example}

\newtheorem{remark}[theorem]{Remark}

\input{epsf.sty}
\usepackage{graphicx}

\begin{document}

\title{Connected Quandles Associated with Pointed Abelian Groups}

\author{W. Edwin Clark, Mohamed Elhamdadi, 
Xiang-dong Hou, \\
Masahico Saito,  Timothy Yeatman\\[2mm]
Department of Mathematics and Statistics\\
University of South Florida}

\date{\empty}

\maketitle

\begin{abstract}
A quandle is a self-distributive  algebraic structure that appears
 in quasi-group and knot theories. 
For each abelian group $A$ and $c \in A$ we  define a quandle $G(A,c)$ on $\Z_3 \times A$.  
These quandles are generalizations of a class of non-medial Latin quandles 
defined by V. M. Galkin so we call them {\it Galkin quandles}. Each $G(A,c)$
is connected but not Latin unless $A$ has odd order. $G(A,c)$ is non-medial unless $3A = 0$.  
We classify their isomorphism classes in terms of  pointed abelian groups, 
and study their various properties. 
A family of symmetric connected quandles is 
constructed from Galkin quandles, and
some aspects of knot colorings by Galkin quandles are also discussed. 
\end{abstract}

\section{Introduction}\label{Introsec}

Sets with certain self-distributive operations called  {\it quandles}
have been studied since 1940s (for example, \cite{Taka})
in various areas with different names. 
They have been studied, for example, 
 as an algebraic system for symmetries
and as quasi-groups.
The {\it fundamental quandle}
was defined in a manner similar to the fundamental group \cite{Joyce,Mat}, which made quandles
an important tool in knot theory.
Algebraic homology theories for quandles have been developed \cite{CJKLS,FRS1},
 extensions
of quandles by cocycles have been studied \cite{AG,CENS}, and
applied to various properties of 
knots and knotted surfaces (\cite{CKS}, for example). 

Before algebraic theories of extensions were developed, Galkin~\cite{galkin}
defined a family of quandles 
 that are extensions of 
the $3$-element connected quandle $R_3$, 
and we call them {\it Galkin quandles}. 
Even though the definition of Galkin quandles is a special case of 
a cocycle extension described   in \cite{AG},
they have curious properties  such as the 
 explicit  and simple 
defining formula, close connections to dihedral quandles, and the fact that 
they appear in the list of small connected quandles.

In this paper, we generalize Galkin's definition and define a family of quandles
that are extensions of $R_3$,
characterize their isomorphism classes, and study their properties.
The definition is given in Section~\ref{defsec}
after a brief review of necessary materials in Section~\ref{prelimsec}.
Isomorphism classes are characterized by pointed abelian groups in Section~\ref{isomsec}.
Various algebraic properties of Galkin quandles are investigated 
in Section~\ref{propertysec},
and their knot colorings are studied in Section~\ref{knotsec}.

\smallskip

{\it Acknowledgement:}
Special thanks to  Michael Kinyon for bringing Galkin's paper~\cite{galkin}
to our attention and pointing out  the
 construction of non-medial, 
 Latin quandles on page $950$
of \cite{galkin}
 that we call here Galkin quandles $G(\Z_p,c_1,c_2)$.  
We are also grateful to Professor Kinyon for helping us with using Mace4
for colorings of knots by quandles,
and for telling us about Belousov's work on distributive quasigroups. 
Thanks to David Stanovsky for useful discussions on these matters.
We are  grateful to James McCarron for his help with the Magma package in {\it Maple 15} especially with isomorphism testing. 
M.S. was supported in part by NSF grant DMS  \#0900671.

\section{Preliminaries}\label{prelimsec}

In this section we briefly review some definitions and examples of quandles. 
More details can be found, for example, in \cite{AG,CKS,FRS1}. 

A {\it quandle} $X$ is a non-empty set with a binary operation $(a, b) \mapsto a * b$
satisfying the following conditions.
\begin{eqnarray}
\mbox{\rm (Idempotency) } & &  \mbox{\rm  For any $a \in X$,
$a* a =a$.} \label{axiom1} \\
\mbox{\rm (Invertibility)}& & \mbox{\rm For any $b,c \in X$, there is a unique $a \in X$ such that 
$ a*b=c$.} \label{axiom2} \\
\mbox{\rm (Right self-distributivity)} & & 
\mbox{\rm For any $a,b,c \in X$, we have
$ (a*b)*c=(a*c)*(b*c). $} \label{axiom3} 
\end{eqnarray}
 A {\it quandle homomorphism} between two quandles $X, Y$ is
 a map $f: X \rightarrow Y$ such that $f(x*_X y)=f(x) *_Y f(y) $, where
 $*_X$ and $*_Y$ 
 denote 
 the quandle operations of $X$ and $Y$, respectively.
 A {\it quandle isomorphism} is a bijective quandle homomorphism, and 
 two quandles are {\it isomorphic} if there is a quandle isomorphism 
 between them.
 
 Typical examples of quandles include the following. 
  \begin{itemize}
 \setlength{\itemsep}{-3pt}
\item
Any non-empty set $X$ with the operation $x*y=x$ for any $x,y \in X$ is
a quandle called the {\it trivial} quandle.

\item
A group $X=G$ with
$n$-fold conjugation
as the quandle operation: $a*b=b^{-n} a b^n$.

\item
Let $n$ be a positive integer.
For  
$a, b \in \Z_n$ (integers modulo $n$), 
define
$a\ast b \equiv 2b-a \pmod{n}$.
Then $\ast$ defines a quandle
structure  called the {\it dihedral quandle},
  $R_n$.
This set can be identified with  the
set of reflections of a regular $n$-gon
  with conjugation
as the quandle operation.
\item
Any ${\Z }[T, T^{-1}]$-module $M$
is a quandle with
$a*b=Ta+(1-T)b$, $a,b \in M$, called an {\it  Alexander  quandle}.
An Alexander quandle is also regarded as a pair $(M, T)$ where 
$M$ is an abelian group and $T \in {\rm Aut}(M)$.
  \end{itemize}
  
  Let $X$ be a quandle.
  The {\it right translation}  ${\cal R}_a:X \rightarrow  X$, by $a \in X$, is defined
by ${\cal R}_a(x) = x*a$ for $x \in X$. Similarly the {\it left translation} ${\cal L}_a$
is defined by ${\cal L}_a(x) = a*x$. Then ${\cal R}_a$ is a permutation of $X$ by Axiom (2).
The subgroup of ${\rm Sym}(X)$ generated by the permutations ${\cal R}_a$, $a \in X$, is 
called the {\it inner automorphism group} of $X$,  and is 
denoted by ${\rm Inn}(X)$. 
 We list some definitions of commonly known properties of quandles below.

\begin{itemize}
\setlength{\itemsep}{-3pt}
\item
A quandle is {\it connected} if ${\rm Inn}(X)$ acts transitively on $X$.

\item
A {\it Latin quandle} is a quandle such that for each $a \in X$, the left translation ${\cal L}_a$ is a bijection. 
That is, the multiplication table of the quandle is a Latin square.

\item
A quandle is {\it faithful} if the mapping $a \mapsto {\cal R}_a$ is 
an injection 
from $X$ to ${\rm Inn}(X)$.

\item
A quandle $X$ is {\it involutory}, or a {\it kei}, 
if the right translations  are  involutions: ${\cal R}_a^2 ={\rm  id}$ for all $a \in X$.

\item
It is seen that the operation $\bar{*}$ on $X$ defined by $a\ \bar{*}\ b= {\cal R}_b^{-1} (a) $
is a quandle operation, and $(X,  \bar{*}) $ is called the {\it dual} quandle of $(X, *)$.
If  $(X,  \bar{*}) $ is isomorphic to  $(X, *)$, then $(X, *)$ is called {\it self-dual}.

\item
A quandle $X$ is {\it medial} if 
$(a*b)*(c*d)=(a*c)*(b*d)$  for all  $a,b,c,d \in X$. 
It is also called {\it abelian}. 
It is known and easily seen that 
every Alexander quandle is medial.

\end{itemize}

\begin{figure}[ht]
\begin{minipage}[b]{0.5\linewidth}
\centering
\includegraphics[scale=.5]{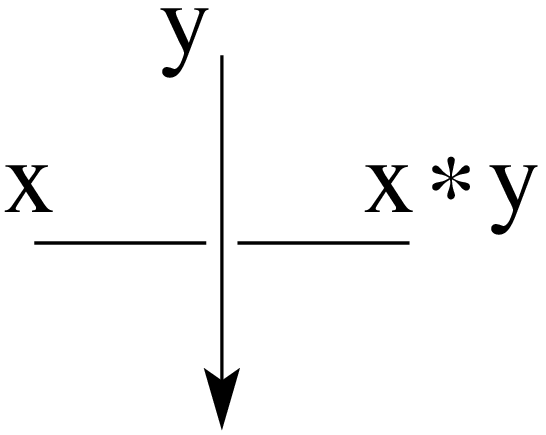}
\caption{A coloring rule at a crossing}
\label{color}
\end{minipage}
\begin{minipage}[b]{0.5\linewidth}
\centering
\includegraphics[scale=.4]{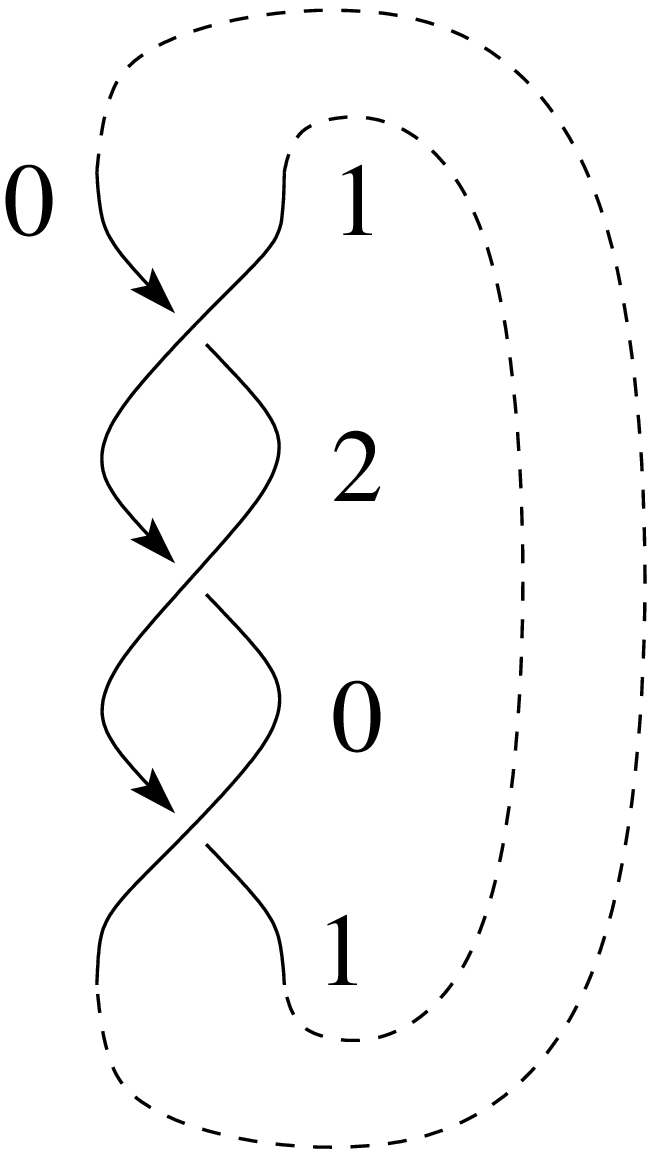}
\caption{Tefoil as the closure of $\sigma_1^3$}
\label{tref}
\end{minipage}
\end{figure}

A {\it coloring} of 
an oriented 
 knot diagram by a quandle $X$ is a map
${\cal C}: {\cal A} \rightarrow X$ 
 from the set of arcs 
 ${\cal A}$ 
of the diagram to $X$ such that the image of the map satisfies the relation depicted in Figure~\ref{color} at each crossing.
More details can be found in \cite{CKS}, for example. 
A coloring that assigns the same element of $X$ for all the arcs is called trivial,
otherwise non-trivial. 
The number of colorings of a knot diagram by a finite quandle is known to be independent of
the choice of a diagram, and hence is a knot invariant. 
A coloring by a dihedral quandle $R_n$ for a positive integer $n>1$ is called
an $n$-coloring. 
If a knot is non-trivially colored by a dihedral quandle $R_n$ for a positive integer $n>1$, 
then it is called {\it $n$-colorable}. 
In Figure~\ref{tref}, a non-trivial $3$-coloring of the trefoil knot ($3_1$ in a common notation 
in a knot table \cite{KI}) is indicated.
This is presented in a closed braid form. 
Each crossing corresponds to a standard generator $\sigma_1$ 
of 
the $2$-strand braid group, and 
$\sigma_1^3$ represents three crossings together  as in the figure. 
The dotted line indicates the closure, see \cite{Rolf} for more details of braids.

The fundamental quandle
is defined in a manner similar to the fundamental group \cite{Joyce,Mat}.
A {\it presentation} of a quandle is defined in a manner similar to
groups as well, and a presentation of the fundamental quandle 
is obtained from  a knot diagram (see, for example, \cite{FR}),
by assigning generators to arcs of a knot diagram, and relations corresponding to crossings. 
The set of a coloring of a knot diagram $K$ by a quandle $X$, then, is in one-to-one 
correspondence with the set of quandle homomorphisms 
from the fundamental quandle of $K$ to $X$.

\section{Definition and notation for Galkin quandles } \label{defsec}

Let  $A$ be an abelian group, also regarded naturally as a $\Z$-module. 
Let $\mu: \Z_3  \rightarrow \Z$ ,  $\tau: \Z_3 \rightarrow A$ be functions. 
Define a binary operation on $\Z_3 \times A$ by 
$$(x, a)*(y, b)=(2y-x, -a + \mu(x-y) b + \tau(x-y) ) \quad x, y \in \Z_3, \ a, b \in A. $$

\begin{proposition} \label{galprop}
 For any abelian group $A$, 
the above operation $*$ defines a quandle structure on $\Z_3 \times A$ if 
$\mu(0)=2$, $\mu(1)=\mu(2)=-1$, and $\tau(0)=0$. 
\end{proposition}

Galkin gave this definition in \cite{galkin}, page $950$, for $A = \Z_p$.
The proposition generalizes his result to any abelian group $A$.
For the proof, 
we examine the axioms. 

\begin{lemma}\label{idemlem}
$(A)$ The operation is idempotent $($i.e., satisfies Axiom~\eqref{axiom1}$)$ if and only if $\mu(0)=2$ and $\tau(0)=0$. 
 $(B)$  The operation as right action is invertible $($i.e. satisfies Axiom~\eqref{axiom2}$)$. 
\end{lemma}
{\it Proof.\/} Direct calculations. $\Box$

\begin{lemma}\label{mutaulem}
The operation $*$ on $\Z_3 \times A$   
is right self-distributive $($i.e., satisfies Axiom~\eqref{axiom3}$)$
if and only if $\mu, \tau$ satisfy the following conditions for any $X, Y \in \Z_3$:
\begin{eqnarray}
\mu(-X)&=& \mu (X) , \label{mu1} \\
\mu (X+Y) + \mu(X - Y)&=&  \mu(X) \mu(Y) , \label{mu2} \\
\tau(X+Y) + \tau(Y-X) &=&  \tau(X) + \tau(-X) + \mu(X) \tau(Y) . \label{mu3}
\end{eqnarray}
\end{lemma}
{\it Proof.\/}
The right self-distributivity
$$(( x, a)*(y,b))*(z,c) = ((x,a)*(z,c))*((y,b)*(z,c))\quad {\rm  for } \quad x,y,z \in \Z_3
\quad {\rm  and } \quad
a,b,c \in A, $$ 
is satisfied if and only if 
\begin{eqnarray*} 
\mu (x-y) &=& \mu(y-x) , \\
\mu (2y - x - z) &=& - \mu (x-z) + \mu(y-x) \mu(y-z), \\
- \tau(x-y) + \tau(2y-x-z) &=& - \tau(x-z) + \mu(y-x) \tau(y-z) + \tau(y-x),
\end{eqnarray*}
by equating the coefficients of $b$, $c$, and the constant term.

The first one is equivalent to \eqref{mu1} by setting $X=x-y$. 
The second is equivalent to \eqref{mu2}  by setting
$X=y-x$ and $Y=z-y$. The third is equivalent to \eqref{mu3} by 
$X=y-x$ and $Y=y-z$. $\Box$

\bigskip

\noindent
{\it Proof of Proposition~\ref{galprop}.\/}
Assume the conditions stated.
By Lemma~\ref{idemlem}
Axioms~\eqref{axiom1} and \eqref{axiom2} are satisfied 
under the specifications 
$\mu(0)=2, \mu(1)=\mu(2)=-1$, and $\tau(0)=0$.

If $X=0$ or $Y=0$, then \eqref{mu2} (together with \eqref{mu1}) 
becomes tautology. 
If $X-Y=0$ or $X+Y=0$, then \eqref{mu2} reduces to $\mu(2X)+2=\mu(X)^2$ 
which is satisfied by the above specifications. 
For $R_3$, if $X+Y\neq 0$ and $X - Y \neq 0$, then either $X=0$ or $Y=0$.
Hence \eqref{mu2} is satisfied.
For \eqref{mu3}, it is checked similarly, for the two cases [$X=0$ or $Y=0$],
and [$X-Y=0$ or $X+Y=0$]. 
$\Box$

\begin{definition} 
{\rm
Let $A$ be an abelian group. 
The quandle defined by $*$ on $\Z_3 \times A$ by Proposition~\ref{galprop}
with 
$\mu(0)=2, \mu(1)=\mu(2)=-1$ and $\tau(0)=0$ 
is called the {\it Galkin quandle  }
and denoted by $G(A, \tau)$. 

Since $\tau$ is specified by the values $\tau(1)=c_1$ and $\tau(2)=c_2$, 
where $c_1, c_2 \in A$, 
we also denote it by $G(A, c_1, c_2)$. 
}
\end{definition}

\begin{lemma}\label{c1c2lem}
For any abelian group $A$  and $c_1, c_2 \in A$, 
$G(A,  c_1, c_2)$ and $ G(A,  0, c_2 - c_1 ) $
are isomorphic. 
\end{lemma}
{\it Proof.\/}
 Let $c = c_2 - c_1$. 
Define 
$\eta: G(A, c_1, c_2) \rightarrow G (A,  0, c)$,
as a map on $\Z_3 \times A $, 
by 
$\eta(x,a)= ( x, a+ \beta(x) )$,
where $\beta(0)=\beta(1)=0$, and $\beta(2)=-c_1$. 
This $\eta$ is a bijection, and we show that it is a quandle homomorphism. 
We compute $\eta( (x, a)*(y,b) )$ and $\eta(x,a) * \eta(y,b)$ for $x,y \in \Z_3$, $a,b \in A $.

If  $x=y$, then  $\mu(x-y)=2$ and $\tau(x-y)=0$ for both 
$G (A,  c_1,  c_2)$ and  $G (A,   0, c)$, so that 
\begin{eqnarray*}
\eta( (x, a)*(x,b) )&=&\eta( x, 2b-a)=(x, 2b-a + \beta(x) ), \\ 
\eta(x,a) * \eta(x,b) &=& (x, a + \beta(x) )*( x, b + \beta(x) )
=(\ x, \ 2(  b  +  \beta(x) ) - ( a +  \beta(x) ) \ ) \\
 &=& ( x, 2b-a  +  \beta(x) ) 
\end{eqnarray*}
as desired.

If $x-y=1 \in \Z_3$, then $\mu(x-y)=-1$ for  both 
$G (A,  c_1,  c_2)$ and  $G (A,  0, c)$
and $\tau(x-y)= c_1$ for $G(A,  c_1,  c_2)$
but $\tau(x-y)=0$ for $G (A,   0, c)$, so that
\begin{eqnarray*}
\eta( (x, a)*(y,b) ) &=& \eta( 2y - x ,-a-b +  c_1 )=(2y-x, -a-b +  c_1 + \beta(2y-x), \\
\eta(x,a) * \eta(y, b) &=& (x, a+ \beta(x) )*( y, b + \beta(y)  )
=(\ 2y-x, \ - ( a + \beta(x) )- ( b+ \beta(y)  ) \ ) . 
\end{eqnarray*}
This holds if and only if $\beta(x) + \beta(y) + \beta(2y-x)=- c_1$, 
which is true since $x\neq y$ implies that exactly one of 
$x, y, 2y-x$ is $2 \in \Z_3$. 

If $x-y=2 \in \Z_3$, then $\mu(x-y)=-1$  for  both 
$G( A, c_1,  c_2)$ and  $G ( A, 0, c)$
and $\tau(x-y)= c_2$ for $G ( A, c_1,  c_2)$
but $\tau(x-y)=c_2 - c_1=c$ for $G ( A, 0, c)$, so that
\begin{eqnarray*}
\eta( (x, a)*(y,b) )&=&\eta( 2y - x ,-a-b +  c_2 )=(2y-x, -a-b +  c_2 + \beta(2y-x) ), \\
\eta(x,a) * \eta(y, b) &=& (x, a + \beta(x)  )*( y, b+ \beta(y)  )
=(2y-x, - ( a+ \beta(x) )- ( b+ \beta(y) ) ) \\
&=& (2y-x, -a-b  - \beta(x) - \beta(y) + (c_2 - c_1)  ),
\end{eqnarray*}
and again 
these are equal 
for the same reason as above.
$\Box$

\bigskip

\noindent
{\bf Notation.}
Since by 
 Lemma~\ref{c1c2lem}, any Galkin quandle is isomorphic to 
$G(A, 0, c)$ for an abelian group $A$ and $c \in A$, we denote
$G(A, 0, c)$ by $G(A, c)$ for short. 

Any finite abelian  group is a product 
$\Bbb Z_{n_1}\times\cdots\times\Bbb Z_{n_k}$ 
 where the positive integers $n_j$ 
satisfy that $n_j | n_{j+1}$ for $j=1, \ldots, k-1$. 
In this case any element $c \in A$ is written in a vector form
$[c_1, \ldots, c_k]$, where $c_j \in \Z_{n_j}$.
Then the corresponding Galkin quandle is denoted by 
$G(A, [c_1, \ldots, c_k])$. 

\bigskip

 \begin{remark}  \label{funclem} 
{\rm
We note that the definition of Galkin quandles induces a functor.
Let ${\bf Ab}_0$ denote the category of pointed abelian groups, that, is the category whose
objects are pairs $(A,c)$ where $A$ is an abelian group and $c \in A$ and whose morphisms
$f: (A,c) \rightarrow (B,d)$ are group homomorphisms $f:A\rightarrow B$ such that $f(c) = d$. 
Let ${\bf Q}$ be the category of quandles, consisting quandles as objects and 
quandle homomorphisms as morphisms.

Then the  correspondence $ (A, c) \stackrel{ {\cal F} }{\mapsto} G(A, c)$ defines a functor 
${\cal F} : {\bf  Ab }_0 \rightarrow {\bf  Q} $. 
It is easy to verify that if  a morphism $f:(A,c)\rightarrow (B,d)$ is given
then the mapping ${\cal F}(f)(x,a) = (x,f(a)) $, $(x,a) \in G(A,c)=\Z_3 \times A$,  is a homomophism from
$G(A,c)$ to $G(B,d)$, and  satisfies ${\cal F}(gf) = {\cal F}(g) {\cal F}(f)$  and 
${\cal F} (\id_{(A,c)} ) = \id_{G(A,c)}$. 
}
\end{remark}

\section{Isomorphism classes} \label{isomsec}

In this section we classify isomorphism classes of Galkin quandles.

\begin{lemma}\label{isomlem}
Let $A$ be an abelian group, and $h: A \rightarrow A'$ be 
a group  isomorphism. 
Then Galkin quandles $G(A, \tau)$ and $G(A' , h \tau)$ are isomorphic as quandles.
\end{lemma}
{\it Proof.\/}
Define 
$f: G(A, \tau)\rightarrow G(A' , h\tau )$,
as a map from 
$\Z_3 \times A$
to  $\Z_3 \times A'$, 
by 
$f(x,a)= ( x, h(a) )$. 
This $f$ is a bijection, and we show that it is a quandle homomorphism by 
computing  $f( (x, a)*(y,b) )$ and $f(x,a) * f(y,b)$ for $x,y \in \Z_3$, $a,b \in A$.

If  $x=y$, then  $\mu(x-y)=2$ and $\tau(x-y)=0=h\tau(0)$ for both 
$G( A, \tau)$ and  $G ( A' , h \tau)$, so that 
\begin{eqnarray*}
f( (x, a)*(x,b) )&=&f( x, 2b-a)=(x, h(2b-a ) ), \\ 
f(x,a) * f(x,b) &=& (x,  h(a) )*( x, h(b) )
=(\ x, \ 2 h(b)-h(a)   \ )
\end{eqnarray*}
as desired.

If $x-y=1 \in \Z_3$, then $\mu(x-y)=-1$ for  both 
$G ( A, \tau)$ and  $G ( A' , h \tau )$.
\begin{eqnarray*}
f( (x, a)*(y,b) ) &=& f( 2y - x ,-a-b + \tau(x-y) )=(\ 2y-x,\  h(-a-b+ \tau(x-y) ) \ )  , \\
f(x,a) * f(y, b) &=& (x, h(a) )*( y, h(b)  )
=(\ 2y-x, \ - h(a) - h(b) + h\tau (x-y) \ )
\end{eqnarray*}
as desired. 
$\Box$

\begin{lemma}\label{ctodlem}
Let $c, d, n$ be positive integers. 
If $ {\rm gcd}(c,n) = d$,  then $G(\Z_n ,  c)$ is isomorphic to $G(\Z_n, d )$.
\end{lemma}
{\it Proof.\/}
 If $A = \Z_n$  then ${\rm Aut}(A) = \Z_n^* = {\rm units\  of }$ $ \Z_n$,  
 and the divisors of $n$
are representatives of the orbits of  $\Z_n^*$ acting on $\Z_n$. 
$\Box$

\bigskip

Thus  we may choose the divisors of $n$  for
the values of $c$ for representing isomorphism classes of $G(\Z_n, c)$. 

\begin{corollary}\label{pgpcor}
 If $A$  is a vector space $($elementary $p$-group$)$ then there are 
 exactly two isomorphism classes 
 of Galkin quandles $G(A, \tau)$. 
\end{corollary}
{\it Proof.\/}
 If $A$ is a vector space containing non-zero vectors $c_1$ and $c_2$, 
 then 
   there is
a non-singular linear transformation $h$ of $A$ such that $h(c_1) = c_2$. 
That $G(A,0)$ is not isomorphic to $G(A,c)$ if $c \neq 0$ follows from
Lemma~\ref{cyclelem}  below.
$\Box$

\bigskip
 
For distinguishing isomorphism classes, cycle structures of the right action 
is useful, and we use the following lemmas.

\begin{lemma}\label{connlem}
For any abelian group $A$,   the Galkin quandle
 $G(A, \tau)$ 
 is connected. 
\end{lemma}
{\it Proof.\/}  Recall that the operation  is defined by 
 the formula 
  $$ (x,a)*(y,b) = (2y - x, -a + \mu(x-y)b + \tau(x-y)), $$
$ \mu(0) = 2$, $\mu(1)=\mu(2)=-1$ and $\tau(0)=0$.
If $x\neq y$, then  $ (x,a)*(y,b) = ( 2y-x, -a -b +c_i) = (z, c)$ 
where $i=1$ or $2$ and  $x, y \in \Z_3$, $a,b \in A$. 
Note that $\{ x, y, 2y-x \}=\Z_3$ if $x \neq y$. 
In particular, for any $(x,a)$ and $(z, c)$ with  $x \neq z$, there is 
$(y, b)$ such that  $(x,a) * (y,b)=(z, c)$. 

For any  $(x, a_1)$ and $(x, a_2)$ where $x \in \Z_3$, $a_1, a_2 \in A$, 
take $(z, c) \in \Z_3 \times A$ such that 
$z \neq x$. Then there are $(y, b_1), (y, b_2)$ such that 
$x \neq y \neq z$ and 
$(x, a_1)*(y, b_1)=(z, c)$ and $(z, c)*(y, b_2)=(x, a_2)$. 
Hence $G(A, \tau)$ is connected.
$\Box$

\begin{lemma}\label{cyclelem}
The cycle structure of
a right  translation
 in $G(A,\tau)$  where $\tau(0) = \tau(1) = 0$  and $\tau(2) = c$,
consists of $1$-cycles, 
$2$-cycles and $2k$-cycles where $k$ is the order of $c$ in the group $A$.

Since isomorphic quandles have the same cycle structure of right translations, 
$G(A, c)$ and $G(A, c')$ for $c, c' \in A$ are not isomorphic 
unless  the orders of $c$ and $c'$ coincide.
\end{lemma}
{\it Proof.\/}
Let $\tau(0)=0$, $ \tau(1) = 0$ and $\tau(2) = c$. Then by Lemma~\ref{connlem}, 
 the cycle
structure of each column  is the same as the cycle
structure of the right translation by $(0,0)$,  that is, of the
permutation $f(x,a) = (x,a)*(0,0) = (-x,-a+\tau(x))$.

We show that this permutation  has cycles of length only $1$, 
$2$ and twice 
the order of $c$ in $A$. 
Since $f(0,a) = (0,-a) $ for $a \in A$, $a \neq 0$, 
we have $f^2 (0, a)=(0,a)$ so that $(0,a)$ generates a $2$-cycle,
or a $1$-cycle if $2a = 0$. 
Now
from $f(1,a) = (2,-a )$ and 
$f(2,a) = (1,-a + c)$ for $a \in A$, 
by induction it is easy to see that for $k > 0$, 
$f^{2k} (1,a) = (1,a + kc)$ and 
$f^{2k} (2,a) = (2,a - kc)$. 
In the case of $(1,a)$, $a \neq 0$,  the cycle closes when $a + kc = a$ in $A$.
The smallest $k$ for which this holds is the order of $c$, 
in which case the cycle is of length $2k $.
A cycle beginning at $(2,a)$ similarly has this same
length. 
$\Box$

\begin{proposition}\label{cycisomprop}
Let $n$ be a positive integer, $A=\Z_n$, and $c_i, c_i' \in \Z_n$
for $i=1,2$. 
Two Galkin quandles $G(A, c_1,  c_2)$ and $ G(A, c_1', c_2') $
 are isomorphic if and only if 
$\gcd (c_1 - c_2, n)=\gcd (c_1' - c_2', n)$.
\end{proposition}
{\it Proof.\/}
If $\gcd (c_1 - c_2, n)=\gcd (c_1' - c_2', n)$, 
then they are isomorphic by Lemmas~\ref{c1c2lem} and \ref{ctodlem}. 
The cycle structures are different if 
$\gcd (c_1 - c_2, n)\neq \gcd (c_1' - c_2', n)$
by Lemma~\ref{cyclelem},  and hence they are not isomorphic.
$\Box$

\begin{remark}{\rm 
The cycle structure is not sufficient for non-cyclic groups $A$.
For example, let $A=\Z_2 \times \Z_4$. 
 Then $G(A, [1,0])$ and $G(A, [0,2])$
 have the same cycle structure 
 for right translations,  
 with cycle lengths $\{2,2,4,4,4,4\}$ in a multiset notation,
 yet they are known to be not isomorphic.
 (In the notation of Example~\ref{rigex} (see below), 
 $G(A, [1,0])=C[24, 29]$ and $G(A, [0,2])=C[24, 31]$, that are not isomorphic.)
 We note that there is no automorphism of $A$ carrying $[1,0]$ to $[0,2]$. 

}
\end{remark}

More generally the isomorphism classes of Galkin quandles are characterized as follows.

\begin{theorem}\label{isomthm}
{\rm
Suppose $A, A'$ are finite abelian groups. 
Two Galkin quandles 
$G(A,  \tau)$ and $G(A', \tau')$ are isomorphic if and only if there exists 
a group  isomorphism $h: A \rightarrow A'$ such that $h\tau=\tau'$. 
}
\end{theorem}

One implication is Lemma~\ref{isomlem}. For the other, 
first we prove the following two lemmas.
We will use a well known description of the automorphisms of a finite abelian group which can be found in \cite{ Hil-Rhe07, Ran1907}.

\begin{lemma}\label{L1}
Let $A$ be a finite abelian $p$-group and let $f:pA\to pA$ be an automorphism. Then $f$ can be extended to an automorphism of $A$.
\end{lemma}
{\it Proof.\/} 
Let $A=\Bbb Z_{p^1}^{n_1}\times \cdots\times\Bbb Z_{p^k}^{n_k}$. Then 
\begin{equation}\label{1}
f(\left[\begin{matrix} px_2\cr \vdots\cr px_k\end{matrix}\right])=P\left[\begin{matrix} px_2\cr \vdots\cr px_k\end{matrix}\right],
\qquad \left[\begin{matrix} x_2\cr \vdots\cr x_k\end{matrix}\right]\in \Bbb Z_{p^2}^{n_2}\times \cdots\times\Bbb Z_{p^k}^{n_k},
\end{equation}
where
\begin{equation}\label{2}
P=\left[\begin{matrix} P_{22}&P_{23}&\cdots&P_{2k}\cr
pP_{32}&P_{33}&\cdots&P_{3k}\cr
\vdots&\vdots&&\vdots\cr
p^{k-2}P_{k2}&p^{k-3}P_{k3}&\cdots&P_{kk}\end{matrix}\right],
\end{equation}
$P_{ij}\in\text{M}_{n_i\times n_j}(\Bbb Z)$, $\det P_{ii}\not\equiv 0\pmod p$.
Define $g:A\to A$ by
\[
g(\left[\begin{matrix} x_1\cr x_2\cr \vdots\cr x_k\end{matrix}\right])=\left[\begin{matrix}I\cr &P\end{matrix}\right]\left[\begin{matrix} x_1\cr x_2\cr \vdots\cr x_k\end{matrix}\right],
\qquad \left[\begin{matrix} x_1\cr x_2\cr \vdots\cr x_k\end{matrix}\right]\in \Bbb Z_{p^1}^{n_1}\times\Bbb Z_{p^2}^{n_2}\times \cdots\times\Bbb Z_{p^k}^{n_k}.
\]
Then $g\in\text{Aut}(A)$ and $g|_{pA}=f$.
$\Box$ 

\begin{lemma}\label{L2}
Let $A$ be a finite abelian $p$-group and let $a,b\in A\setminus pA$. If there exists an automorphism $f:pA\to pA$ such that $f(pa)=pb$, then there exists an automorphism $g:A\to A$ such that $g(a)=b$.
\end{lemma}
{\it Proof.\/} 
Let $A=\Bbb Z_{p^1}^{n_1}\times \cdots\times\Bbb Z_{p^k}^{n_k}$ and let $f$ be defined by \eqref{1} and \eqref{2}. Write
\[
a=\left[\begin{matrix} a_1\cr \vdots\cr a_n\end{matrix}\right],\ b=\left[\begin{matrix} b_1\cr \vdots\cr b_n\end{matrix}\right],\quad a_i,b_i\in\Bbb Z_{p^i}^{n_i}.
\]
Since $f(pa)=pb$, we have
\[
p\Bigl(P \left[\begin{matrix} a_2\cr \vdots\cr a_n\end{matrix}\right]-\left[\begin{matrix} b_2\cr \vdots\cr b_n\end{matrix}\right]\Bigr)=0,
\]
i.e.,
\begin{equation}\label{3}
P \left[\begin{matrix} a_2\cr \vdots\cr a_n\end{matrix}\right]-\left[\begin{matrix} b_2\cr \vdots\cr b_n\end{matrix}\right]=
\left[\begin{matrix} pc_2\cr \vdots\cr p^{k-1}c_k\end{matrix}\right],\quad c_i\in \Bbb Z_{p^i}^{n_i},\ 2\le i\le k.
\end{equation}

\medskip

{\bf Case 1.} Assume that 
$
\left[\begin{matrix} a_2\cr \vdots\cr a_n\end{matrix}\right]\in pA.
$ 
Then by \eqref{3},
$ 
\left[\begin{matrix} b_2\cr \vdots\cr b_n\end{matrix}\right]\in pA.
$ 
So $a_1\ne 0$ and $b_1\ne 0$. Then we have
\[
\left[\begin{matrix} pc_2\cr \vdots\cr p^{k-1}c_k\end{matrix}\right]=\left[\begin{matrix} pQ_2\cr \vdots\cr p^{k-1}Q_k\end{matrix}\right]a_1
\]
for some $Q_i\in\text{M}_{n_i\times n_1}(\Bbb Z)$, $2\le i\le k$. Also, there exists $P_{11}\in\text{M}_{n_1\times n_1}(\Bbb Z)$ such that $\det P_{11}\not\equiv 0\pmod p$
and $P_{11}a_1=b_1$. 
Let $g\in\text{Aut}(A)$ be defined by
\[
g(\left[\begin{matrix} x_1\cr x_2\cr \vdots\cr x_k\end{matrix}\right])=
\left[\begin{matrix} P_{11}&0\cr
\begin{matrix} -pQ_2\cr \vdots\cr -p^{k-1}Q_k\end{matrix} & \kern 3mm P \kern 3mm \end{matrix}\right]
\left[\begin{matrix} x_1\cr x_2\cr \vdots\cr x_k\end{matrix}\right],\qquad x_i\in\Bbb Z_{p^i}^{n_i}.
\]
Then $g(a)=b$.

\medskip

{\bf Case 2.} Assume that 
$ 
\left[\begin{matrix} a_2\cr \vdots\cr a_n\end{matrix}\right]\notin pA.
$ 
Then there exists $2\le s\le k$ such that $a_s\notin p\Bbb Z_{p^s}^{n_s}$. Then we have
 \[
\left[\begin{matrix} c_2\cr \vdots\cr p^{k-2}c_k\end{matrix}\right]=\left[\begin{matrix} Q_2\cr \vdots\cr p^{k-2}Q_k\end{matrix}\right]a_s
 \]
for some $Q_i\in\text{M}_{n_i\times n_s}(\Bbb Z)$, $2\le i\le k$.
Put
\[
Q=\left[\begin{matrix} 0&\cdots&0 & Q_2 & 0&\cdots&0\cr
\vdots& &\vdots &\vdots & \vdots& & \vdots\cr
0&\cdots&0 & p^{k-2}Q_k & 0&\cdots&0
\end{matrix}\right],
\]
where the $(i,j)$ block is of size $n_i\times n_j$ and $Q_2$ is in the $(1,s)$ block. Then
$ 
Q\left[\begin{matrix} a_2\cr \vdots\cr a_k\end{matrix}\right]=\left[\begin{matrix} c_2\cr \vdots\cr p^{k-2}c_k\end{matrix}\right].
$ 
Also, there exist $U\in\text{M}_{n_1\times(n_2+\cdots+n_k)}(\Bbb Z)$ such that
$
U\left[\begin{matrix} a_2\cr \vdots\cr a_k\end{matrix}\right]=b_1-a_1.
$ 
Now define $g\in\text{Aut}(A)$ by 
\[
g(\left[\begin{matrix} x_1\cr x_2\cr \vdots\cr x_k\end{matrix}\right])=
\left[\begin{matrix} I&U\cr
0&P-pQ \end{matrix}\right]
\left[\begin{matrix} x_1\cr x_2\cr \vdots\cr x_k\end{matrix}\right],\qquad x_i\in\Bbb Z_{p^i}^{n_i}.
\]
Then $g(a)=b$.
$\Box$ 

\bigskip

\noindent
{\it Proof of Theorem~\ref{isomthm}.\/} 
We assume that $|3A'|\le |3A|$.
Since $G(A',c')$ is connected, there exists an isomorphism $\phi:G(A,c)\to G(A',c')$ such that $\phi(0,0)=(0,0)$. Write
\[
\phi(x,a)=(\alpha(x,a), \beta(x,a)),\quad (x,a)\in\Bbb Z_3\times A.
\] 
Define $t: \Z_3 \rightarrow A$ by 
$$
t(x)=\begin{cases}
1&\text{if}\ x=2,\cr 
0&\text{if}\ x\ne 2, 
\end{cases}
$$
so that for $(x,a),(y,b)\in\Bbb Z_3\times A$, the operation on $G(A, c)$  is written by 
\[
(x,a)*(y,b)=(-x-y, -a+\mu(x-y)b+ t(x-y)c). 
\]
Then $\phi((x,a)*(y,b))=\phi(x,a)*\phi(y,b)$ is equivalent to
\begin{equation}\label{e1}
\alpha(-x-y,-a+\mu(x-y)b+t(x-y)c)=-\alpha(x,a)-\alpha(y,b),
\end{equation}
\begin{equation}\label{e2}
\begin{split}
&\beta(-x-y, -a+\mu(x-y)b+t(x-y)c)\cr
=\;&-\beta(x,a)+\mu(\alpha(x,a)-\alpha(y,b))\beta(y,b)+t(\alpha(x,a)-\alpha(y,b))c'.
\end{split}
\end{equation}

$1^\circ$ 
We claim that $\alpha(0,\cdot):A\to \Bbb Z_3$ is a homomorphism.

\medskip

\noindent
Setting $x=y=0$ in \eqref{e1} we have
\begin{equation}\label{e3}
\alpha(0,-a+2b)=-\alpha(0,a)-\alpha(0,b).
\end{equation}
Setting $b=0$ in \eqref{e3} we have
\begin{equation}\label{e4}
\alpha(0,-a)=-\alpha(0,a).
\end{equation}
By the symmetry of the RHS of \eqref{e3}, we also have
\begin{equation}\label{e5}
\alpha(0,-a+2b)=\alpha(0,-b+2a),\quad a,b\in A.
\end{equation}
Now we have
\[
\begin{split}
\alpha(0,a+b)\;&=\alpha(0,a-b+2b)\cr
&=\alpha(0,-b+2(b-a))\kern 1.3cm\text{(by \eqref{e5})}\cr
&=\alpha(0,b-2a)\cr
&=-\alpha(0,-b)-\alpha(0,-a)\kern 1cm\text{(by \eqref{e3})}\cr
&=\alpha(0,a)+\alpha(0,b)\kern 1.83cm\text{(by \eqref{e4})}.
\end{split}
\]

\medskip

$2^\circ$
We claim that there exists $u\in\Bbb Z_3$ such that
\begin{equation}\label{e6}
\alpha(x,a)=\alpha(0,a)+ux,\quad (x,a)\in\Bbb Z_3\times A.
\end{equation}

\noindent
Setting $x=1$ and $y=0$ in \eqref{e1} we have
\begin{equation}\label{e7}
\alpha(-1,-a-b)=-\alpha(1,a)-\alpha(0,b).
\end{equation}
Setting $b=0$ in \eqref{e7} gives
\begin{equation}\label{e8}
\alpha(-1,-a)=-\alpha(1,a).
\end{equation}
Letting $a=0$ in \eqref{e7} and using \eqref{e8} we get
\begin{equation}\label{A1}
\alpha(1,b)=\alpha(0,b)+\alpha(1,0),\qquad b\in A.
\end{equation}
Equations \eqref{A1} and \eqref{e8} also imply that
\begin{equation}\label{B}
\alpha(-1,-b)=\alpha(0,-b)-\alpha(1,0),\qquad b\in A.
\end{equation}
Let $u=\alpha(1,0)$. Then
\[
\alpha(x,a)=\alpha(0,a)+ux,\qquad (x,a)\in\Bbb Z_3\times A.
\]

\medskip

$3^\circ$
We claim that 
\begin{equation}\label{e9}
\alpha(0,c)=0.
\end{equation}

\noindent
Substituting \eqref{e6} in\eqref{e1} we get
\[
\alpha(0,-a+\mu(x-y)b+t(x-y)c)=-\alpha(0,a)-\alpha(0,b).
\]
Setting $x-y=2$ we have $\alpha(0,c)=0$.

\bigskip

The remaining part of the proof is divided into two cases according as $u$ is zero or nonzero in \eqref{e6}. 

\bigskip

{\bf Case 1.} Assume $u=0$ in \eqref{e6}. 

\medskip

\noindent
We have $\alpha(x,a)=\alpha(0,a)$ for all $(x,a)\in \Bbb Z_3\times A$. We write $\alpha(a)$ for $\alpha(0,a)$. Then \eqref{e2} becomes
\begin{equation}\label{e10}
\beta(-x-y,-a+\mu(x-y)b+t(x-y)c)=-\beta(x,a)+\mu(\alpha(a-b))\beta(y,b)+t(\alpha(a-b))c'.
\end{equation}

\medskip

$1.1^\circ$
We claim that $c=0$.

\medskip

\noindent
Equation \eqref{e10} with $x=1$, $y=0$, $a=b=0$ yields
\[
\beta(-1,0)=-\beta(1,0),
\]
and with $x=-1$, $y=0$, $a=b=0$, it yields
\[
\beta(1,c)=-\beta(-1,0).
\]
Thus $\beta(1,c)=\beta(1,0)$. Since $\alpha(1,c)=0=\alpha(1,0)$, we have $\phi(1,c)=\phi(1,0)$. Thus $c=0$.

\medskip

$1.2^\circ$ We claim that $c'=0$.

\medskip

\noindent
The homomorphism $\alpha: A\to \Bbb Z_3$ must be onto. (Otherwise $\phi$ is not onto.) Choose $d\in A$ such that $\alpha(d)=-1$. Equation~\eqref{e10} with $x=y=0$, $a=d$, $b=0$ gives
\[
\beta(0,-d)=-\beta(0,d)+c',
\]
and with $x=y=0$, $a=-d$, $b=0$, it gives
\[
\beta(0,d)=-\beta(0,-d).
\]
Therefore $c'=0$.

\medskip

$1.3^\circ$ Now \eqref{e10} becomes
\begin{equation}\label{0.16}
\beta(-x-y,-a+\mu(x-y)b)=-\beta(x,a)+\mu(\alpha(a-b))\beta(y,b).
\end{equation}
Setting $y=0$ and $b=0$ in \eqref{0.16} we have
\begin{equation}\label{0.17}
\beta(-x,-a)=-\beta(x,a).
\end{equation}

\medskip

$1.4^\circ$ We claim that
$\beta(0,\cdot):3A\to A'$ is a 1-1 homomorphism.

\medskip

\noindent
Note that $3A\subset\ker\alpha$.
Let $a,b\in 3A$ and $x=-1$, $y=1$ in \eqref{0.16}. We have
\begin{equation}\label{0.18}
\beta(0,-a-b)=-\beta(-1,a)+2\beta(1,b).
\end{equation}
Setting $b=0$ and $a=0$, respectively, in \eqref{0.18} and using \eqref{0.17} we have
\begin{gather}
\label{0.19}
\beta(0,-a)=-\beta(-1,a)+2\beta(1,0)=\beta(1,-a)+2\beta(1,0),\\
\label{0.20}  
\beta(0,-b)=-\beta(-1,0)+2\beta(1,b)=\beta(1,0)+2\beta(1,b).
\end{gather}
Setting $a=b=0$ in \eqref{0.18} we have
\begin{equation}\label{A}
3\beta(1,0)=0.
\end{equation} 
Combining \eqref{0.18} -- \eqref{A} we have
\[
\beta(0, -a-b ) =\beta(0,-a)+\beta(0,-b).
\]
If $a\in 3A$ such that $\beta(0,a)=0$, then $\phi(0,a)=(0,0)$, so $a=0$. Thus $\beta(0,\cdot):3A\to A'$ is 1-1.

\medskip

$1.5^\circ$ We claim that $\beta(0,3b)\in 3A'$ for all $b\in A$.

\medskip

\noindent
Let $x=y=0$ and $a=-b$ in \eqref{0.16}. We have
\[
\begin{split}
\beta(0,3b)\;&=-\beta(0,-b)+\mu(\alpha(-2b))\beta(0,b)\cr
&=\beta(0,b)+\mu(\alpha(b))\beta(0,b)\cr
&\equiv 0 \pmod{3A'} \kern 2.5cm\text{(since $\mu(\alpha(b))\equiv -1\pmod 3$)}.
\end{split}
\]

\medskip

$1.6^\circ$ 
Now $\beta(0,\cdot):3A\to 3A'$ is a 1-1 homomorphism.  Since $|3A'|\le |3A|$, $\beta(0,\cdot):3A\to 3A'$ is an isomorphism. Since $|A|=|A'|$, we have $A\cong A'$. We are done in Case 1.

\bigskip

{\bf Case 2.} Assume $u\ne 0$ in \eqref{e6}.

\medskip

By the proofs of Lemma 3.5 and Proposition 5.5 (see below), 
$(x',a')\mapsto (-x',a'-t(-x')c')$ is an isomorphism from $G(A',c')$ to $G(A',-c')$. Thus we may assume $u=1$ in \eqref{e6}.
We have $\alpha(x,a)=\alpha(0,a)+x$ for all $(x,a)\in \Bbb Z_3\times A$.

\medskip

$2.1^\circ$
We claim that $\beta(0,\cdot):\ker \alpha(0,\cdot)\to A'$ is a 1-1 homomorphism.

\medskip

\noindent
In \eqref{e2} let $a,b\in\ker\alpha(0,\cdot)$ and $x=-1$, $y=1$. We have
\begin{equation}\label{e18}
\beta(0,-a-b)=-\beta(-1,a)-\beta(1,b).
\end{equation}
Eq.uation~\eqref{e18} with $a=-b$ yields
\begin{equation}\label{e19}
\beta(-1,-b)=-\beta(1,b).
\end{equation}
So
\begin{equation}\label{e20}
\beta(0,-a-b)=\beta(1,-a)-\beta(1,b).
\end{equation}
Letting $b=0$ and $a=0$ in \eqref{e20}, respectively, we have
\begin{eqnarray*} 
\beta(0,-a)&=& \beta(1,-a)-\beta(1,0),\\
\beta(0,-b)&=&\beta(1,0)-\beta(1,b).
\end{eqnarray*} 
Thus
\[
\begin{split}
\beta(0,-a)+\beta(0,-b)\;&=\beta(1,-a)-\beta(1,b)\cr
&=\beta(0,-a-b)\kern 1cm \text{(by \eqref{e20})}.
\end{split}
\]
If $a\in\ker\alpha(0,\cdot)$ such that $\beta(0,a)=0$, then $\phi(0,a)=(0,0)$, so $a=0$. Hence $\beta(0,\cdot):\ker \alpha(0,\cdot)\to A'$ is 
1-1.
\medskip

$2.2^\circ$
We claim that $\beta(0,3a)\in 3A'$ for all $a\in A$.

\medskip

\noindent
Setting $x=y=0$ in \eqref{e2} we have
\begin{equation}\label{e23}
\begin{split}
\beta(0,-a+2b)\;&=-\beta(0,a)+\mu(\alpha(0,a-b))\beta(0,b)+t(\alpha(0,a-b))c'\cr
&\equiv -\beta(0,a)-\beta(0,b)+t(\alpha(0,a-b))c' \pmod{3A'}.
\end{split}
\end{equation}
By \eqref{e23},
\[
\beta(0,3a)=\beta(0, 
-a+2(2a))\equiv -\beta(0,a)-\beta(0,2a)+t(\alpha(0,-a))c'\pmod{3A'}
\]
and
\[
\beta(0,2a)=\beta(0,0+2a)\equiv -\beta(0,a)+t(\alpha(0,-a))c'\pmod{3A'}.
\]
Thus $\beta(0,3a)\equiv 0\pmod{3A'}$.

\medskip

$2.3^\circ$
By the argument in $1.6^\circ$, $\beta(0,\cdot):3A\to 3A'$ is an isomorphism and $A\cong A'$.

\medskip

$2.4^\circ$
We claim that $\beta(0,c)=c'$.

\medskip

\noindent
Equation~\eqref{e2} with $x=1$, $y=-1$, $a=b=0$ yields
\[
\begin{split}
\beta(0,c)\;&=-\beta(1,0)-\beta(-1,0)+c'\cr
&=c'\kern 3.5cm \text{(by \eqref{e19})}.
\end{split}
\]

\medskip

$2.4^\circ$
Now we complete the proof in Case 2. Write $A=A_1\oplus A_2$ and $A'=A_1'\oplus A_2'$, where $3\nmid |A_1|$, $3\nmid |A_1'|$, $|A_2|$ and $|A_2'|$ are powers of $3$. Write $c=c_1+c_2$, where $c_1\in A_1$, $c_2\in A_2$. Then $c_1\in A_1\subset\ker\alpha(0,\cdot)$, so $c_2=c-c_1\in\ker\alpha(0,\cdot)$. Since $\beta(0,\cdot):\ker\alpha(0,\cdot)\to A'$ is a homomorphism, we have
\[
c'=\beta(0,c_1)+\beta(0,c_2)=c_1'+c_2',
\]
where $c_1'=\beta(0,c_1)\in A_1'$ and $c_2'=\beta(0,c_2)\in A_2'$. By $2.3^\circ$, $\beta(0,\cdot):A_1\to A_1'$ is an isomorphism. So it suffices to show that there exists an isomorphism $f:A_2\to A_2'$ such that $f(c_2)=c_2'$.

First assume $c_2\in 3A_2$. Then $c_2'\in 3A_2'$. By Lemma~\ref{L1}, the isomorphism $\beta(0,\cdot):3A\to 3A'$ can be extended to an isomorphism $f:A_2\to A_2'$ and we are done.

Now assume that $c_2\in A_2\setminus 3A_2$. We claim that $c_2\in A_2'\setminus 3A_2'$. Assume to the contrary that $c_2'\in 3A_2'$. 
By $2.3^\circ$, there exists $d\in A_2$ such that $\beta(0,3d)=c_2'=\beta(0,c_2)$. By $2.1^\circ$, $c_2=3d$, which is a contradiction.

Note that $\beta(0,\cdot):3A_2\to 3A_2'$ is an isomorphism and
\[
\begin{split}\beta(0,3c_2)\;&=3\beta(0,c_2) \kern1cm \text{(by $2.1^\circ$)}\cr
&=3c_2'.
\end{split}
\]
By Lemma~\ref{L2}, there exists an isomorphism $f:A_2\to A_2'$ such that $f(c_2)=c_2'$.
$\Box$ 

\begin{table}
$$
\begin{array}{lllllllllll}
C[6,1] &  =& G(\Z_2,[0]) & 1084 & &  C[24,28]&= &G(\Z_8,[4]) & 1084   \\
C[6,2] &  =& G(\Z_2,[1]) & 1084 & &  C[24,29]& =&G(\Z_2 \times \Z_4,  [1, 0], [1, 2] )  & 1084 \\ 
C[9,2] &  =& G(\Z_3,[0]) & 1084 & & C[24,30]& =& G(\Z_2 \times \Z_4, [0,0] )  & 1084 \\
C[9,6] &  =& G(\Z_3,[1]) & 1084 & & C[24,31] &  =&G(\Z_2 \times \Z_4, [0,2] )  & 1084  \\ 
C[12,5] &  =& G(\Z_4,[2]) & 1084 & & C[24,32] &=&G(\Z_8,[1]) & 1051 \\
C[12,6] &  =& G(\Z_4,[0]) & 1084 & & C[24,33]& = & G(\Z_2 \times \Z_4, [0, 1], [1, 1] ) & 1051 \\
C[12,7] &  =& G(\Z_4,[1]) & 1051  & & C[24,38] & =&  G(\Z_2 \times \Z_2 \times \Z_2, [0,0,1] ) & 1084 \\
C[12,8] &  =& G(\Z_2\times \Z_2 ,[0,0])  & 1084 & & C[24,39] & =& G(\Z_2 \times \Z_2 \times \Z_2, [0,0,0])  & 1084  \\
C[12,9] &  =& G(\Z_2 \times \Z_2 ,[1,0])   & 1084 & & C[27,2] &  =& G(\Z_3 \times \Z_3 ,[0,0] )  & 1084  \\ 
C[15,5] &  =& G(\Z_5, [1]) & 1440 & &  C[27,12] &  =& G(\Z_9,[3])  & 1084  \\
C[15,6] &  =& G(\Z_5, [0]) & 1512 & & C[27,13] &  =& G(\Z_9,[0])   & 1084  \\ 
C[18,1] &  =& G(\Z_2 \times \Z_3,[0,0]) & 1084 & &  C[27,23] &  =& G( \Z_3 \times \Z_3 ,[1,0])  & 1084  \\
C[18,4] &  =& G(\Z_2 \times \Z_3,[1,0]) & 1084 & &  C[27,55] &  =& G(\Z_9,[1])  & 1084  \\
C[18,5] &  =& G(\Z_2 \times \Z_3,[1,1]) & 1084 & &  C[30,12] &  =& G(\Z_2 \times \Z_5 ,[0,1]) &1440  \\
C[18,8] &  =& G(\Z_2 \times \Z_3,[0,1]) & 1084 & &  C[30,13] &  =& G(\Z_2 \times \Z_5 ,[0,0]) &1512 \\ 
C[21,7] &  =& G(\Z_7,[1]) & 1339 & & C[30,14] &  =& G(\Z_2 \times \Z_5 ,[1,1]) & 1440 \\
C[21,8] &  =& G(\Z_7,[0]) & 1386 & & C[30,15] &  =& G(\Z_2 \times \Z_5 ,[1,0]) & 1512\\ 
C[24,26] &  =& G(\Z_8,[2]) & 1071 & & C[33,10] &  =& G(\Z_{11},[0])  & 1260 \\
C[24,27]& =& G(\Z_8,[0]) &1084 &   &  C[33,11] &  =& G(\Z_{11},[1]) & 1220 \\
\end{array}
$$
\caption{Galkin quandles in the rig table}
\end{table}

\begin{remark}{\rm
The numbers of
 isomorphism classes of 
order $3n$, from $n=1$ to $n=100$, 
are as follows:

\noindent
$
 1, 2, 2, 5, 2, 4, 2, 10, 5, 4, 2, 10, 2, 4, 4, 20, 2, 10, 2, 10,  
   4, 4, 2, 20, 5, 4, 10, 10, 2, 8, 2, 36, 4, 4, 4, 25, 2, 4, 4, 20,  $
   $
  2, 8, 2, 10, 10, 4, 2, 40, 5, 10, 4, 10, 2, 20, 4, 20, 4, 4, 2, 20,  
   2, 4, 10, 65, 4, 8, 2,10, 4, 8, 2, 50, 2, 4, 10, 10, 4, 8, $
   
   \noindent
      $ 2, 40,  
20, 4, 2, 20, 4, 4, 4, 20, 2, 20, 4, 10, 4, 4, 4, 72, 2, 10, 10, 25. 
$
  
In \cite{CH} it is shown that the number $N(n)$ of isomorphism classes of 
Galkin quandles of order $n$  is multiplicative, that is, if ${\rm gcd}(n,m) = 1$ then 
$N(nm) = N(n)N(m)$, so it suffices to find $N(q^n)$ for all prime powers $q^n$. 
In \cite{CH} it is established that $N(q^n) = \sum_{0\le m\le n}p(m)p(n-m)$, 
where $p(m)$ is the number of partitions of the integer $m$. In particular, $N(q^n)$ is independent of the prime $q$. The sequence $n \mapsto N(q^n)$ appears in the 
On-Line Encyclopedia of Integer Sequences (OEIS, \cite{OEIS})
 as sequence $A000712$. 
}
\end{remark}

\begin{example}
\label{rigex}
{\rm 
In \cite{rig}, connected quandles are listed up to order $35$.
For a positive integer $n>1$, let $q(n)$ be the number of isomorphism classes 
of connected quandles of order $n$. 
For a positive integer $n>1$, 
if $q(n)\neq 0$, then we denote by $C[n,i]$ the $i$-th quandle of order $n$ in their list 
 ($1<n\leq 35$, $i= 1, \ldots, q(n)$). 
We note that $q(n)=0$ 
for   $n=2, 14, 22, 26$, and $34$ (for $1< n \leq 35$). 
The quandle  $C[n,i]$ is denoted by $Q_{n,i}$ in \cite{Leandro}
(and they are left-distributive in \cite{Leandro}, so that the matrix of
$C[n,i]$ is the transpose of the matrix of $Q_{n,i}$). 
Isomorphism classes of Galkin quandles
are identified with those in their list  in Table~1.

The $4$-digit numbers to the right of each row  in Table~1 
indicate 
the numbers of 
knots  that are colored non-trivially by these Galkin quandles, 
out of total $2977$ knots in the table \cite{KI} with $12$ crossings or less. 
See Section~\ref{knotsec} for more on this. 

}
\end{example}

\section{Properties of Galkin quandles} \label{propertysec}

In this section, we investigate various properties of Galkin quandles.

\begin{lemma}\label{latinlem}
The Galkin quandle $G (A,  \tau)$ is Latin if and only if 
$|A|$ is odd. 
\end{lemma}
{\it Proof.\/} 
To show that it is Latin if $n$ is odd, first note that $R_3$  is Latin.
Suppose that $(x,a)*(y,b) = (x,a)*(y',b')$. Then we have the equations
\begin{eqnarray}
 & -x + 2y = -x + 2y' \label{latin1} \\
&  -a + \mu(x-y)b + \tau(x-y) = -a + \mu(x-y')b'+\tau(x-y'). \label{latin2}
\end{eqnarray}
From \eqref{latin1}  it follows that $y = y' $ and it follows from
\eqref{latin2}   that $\mu(x-y)b = \mu(x-y)b'$. Now since $|A|$ is odd,  
the left module action of $2$ on $A$ is invertible, 
 and hence $b = b'$. 
If $|A|$ is even there is a non-zero element $b$ of order $2$
and hence $(0,0)*(0,b) = (0,0)*(0,0)$, so the quandle is not Latin.
 $\Box$

\begin{lemma}
Any Galkin quandle  is faithful.
\end{lemma}
{\it Proof.\/}
We show that if  $ (x, a) * (y, b)=(x, a)*(y', b')$ holds for all $(x,a)$, 
then $(y, b)=(y', b')$. 
We have $y=y'$ immediately. 
{}From the second factor 
$$-a + \mu(x-y)b + \tau(x-y)=-a + \mu(x-y)b' + \tau(x-y), $$ 
we have $\mu(x-y)b = \mu(x-y)b'$ for any $x$.
Pick $x$ such that $x \neq y$, then $\mu(x-y)=-1$, hence $b=b'$.
$\Box$

\begin{lemma} \label{subqlem}
If $A' $ is a subgroup of $A$  and $c' $ is in $A'$, 
 then $G(A',c')$  is a subquandle of $G(A,c')$. 
\end{lemma}
{\it Proof.\/} Immediate. $\Box$

\begin{lemma}\label{subRnlem}
 Any Galkin quandle $G(A, \tau)$  consists of three disjoint 
subquandles $\{x\} \times A$ for $x \in \Z_3$, 
and each is a product of dihedral quandles. 
\end{lemma}
{\it Proof.\/} Immediate. $\Box$

\bigskip

\begin{sloppypar}
We  note the following somewhat curious quandles from Lemma~\ref{subRnlem}:
 For a positive integer $k$, $G( \Z_2^k, [0, \ldots, 0])$ 
is a connected quandle that is a disjoint union of three trivial subquandles of order $2^k$.
\end{sloppypar}

\begin{lemma} \label{R3lem}
The Galkin quandle $G(A, \tau)$ has $R_3$ as a subquandle if and only if 
  $\tau=0$ or 
$3$ divides $ |A| $. 
\end{lemma}
{\it Proof.} 
If $A$ is any group and $\tau = 0$,  then 
 $  (x,0)*(y,0)=(2y-x,0)$ for any $x, y \in \Z_3$, so that 
$ \Z_3 \times  \{ 0 \}$  is a subquandle isomorphic to $R_3$.
If $3$ divides $|A|$,  then $A$ has a subgroup $B$  isomorphic to $\Z_3$.
In the subquandle $\{0\} \times  B$,
 $ (0,a)*(0,b) = (0,-a+2b) $ for
 $a, b \in B$, so that 
$ \{0\} \times B$  is a subquandle isomorphic to $R_3$.

Conversely, let $S = \{ (x,a),(y,b),(z,d) \}$  be a subquandle of $G(A,c) $ 
isomorphic  to $R_3$. 
Note that the quandle operation of 
$R_3$ is commutative,  and the product of any  two elements is equal to the third.
We examine two cases.

{\bf Case 1}:  $x = y = z$.
In this case we have 
$$
\begin{array}{lllll}
(x,a)*(x,b) &=&  (x,-a + 2b) &=& (x,d), \\
(x,b)*(x,a) &=&  (x,-b + 2a) &=&  (x,d) . 
\end{array}
$$
Hence we have
$
-a + 2b = - b + 2a$
so that $3(a-b)   =  0 $.
If  there are no elements of order $3$ in $ A$, then we 
 have $a-b = 0$ and so $b = a$. This is a contradiction 
to the fact that $S$  contains $3$ elements, so there is an element of order $3$ 
in $A$, 
hence  
$3$ divides $|A|$. 

{\bf Case 2}: 
$x, y$ and $z$ are all distinct (if two are distinct
then all three are). 
In this case consider $S = \{(0,a),(1,b),(2,d) \}$.
Now we have 
$$
\begin{array}{lllll}
(2,d)*(0,a) &= &(1, -d - a + c) &=&  (1,b),  \\
(0,a)*(2,d) &=& (1, -a - d) &=& (1,b). 
\end{array}
$$
Hence  we have $-d - a + c = -a - d$, 
so that 
$c = 0$, and we have $\tau=0$.
$\Box$

\begin{lemma}\label{leftlem}
 The Galkin quandle $G(A, \tau) $ is  left-distributive
if and only if $3A=0$, i.e., every element of $A$ has order $3$. 
\end{lemma}
{\it Proof.\/}
Let $\tau(1) = c_1, \tau(2) = c_2$. 
Let $a = (0,0)$, $b=(0, \alpha)$ and $c = (1,0)$ for $\alpha \in A$. 
Then we get
$a*(b*c) = (1,\alpha -c_2+c_1)$ and 
$(a*b)*(a*c) = (1,-2\alpha -c_2+c_1)$.
If these are equal, then 
$3\alpha =0$ for any $\alpha \in A$. 

Conversely, suppose that every element of $A$ has order $3$.
Then we have $\mu(x)a = 2a$ for any $a \in A$. Then 
 one computes 
\begin{equation}\label{LHS}
\begin{split}
& (x, a) * [(y, b)*(z,c)]  \cr
& \quad =  (x*(y*z), -a+b+c- \tau(y-z) + \tau (x- y*z) ) , 
\end{split}
\end{equation}
\begin{equation}
\begin{split}
& [ (x, a) *(y,b) ] * [(x,a)*(z,c)] \cr
& \quad =  ( (x*y)*(x*z), -a+b+c -\tau(x-y) - \tau(x-z) + \tau(x*y - x*z) ) .\label{RHS}
\end{split}
\end{equation}
If all $x, y, z$ are distinct, then $x-y=1$ or $x-y=2$, and $x*y=z$, $x*z=y$, $y*z=x$. 
If $x-y=1$, then $z=x+1$ and $y-z=1$, $x-z=2$,
and one computes that \eqref{LHS} $=-c_1=$ \eqref{RHS}.
If $x-y=2$, then one computes \eqref{LHS} $=-c_2=$ \eqref{RHS}.
The other cases for $x,y,z$ are checked similarly. 
$\Box$

\begin{proposition}\label{alexprop}
The  Galkin quandle $G(A, \tau)$ is Alexander if and only if $3A=0$. 
\end{proposition}
{\it Proof.\/}
If  $G(A, \tau)$ is Alexander then
since it is connected  
 it is left-distributive, hence  Lemma~\ref{leftlem}
implies $3A=0$. 
Conversely, suppose $3A=0$. 
Then $A=\Z_3^k$ for some positive integer $k$, and is an elementary $3$-group.
By Corollary~\ref{pgpcor}
there are two isomorphism classes,
$G ( \Z_3^k, [0, \ldots, 0 ])$ and $G(\Z_3^k, [0, \ldots, 0, 1])$. 
The quandle $G(\Z_3, 1)=C[9,6]$ is isomorphic to 
$ \Z_3 [t]/ (t+1)^2$ by a direct comparison.
Hence the two classes are isomorphic to the Alexander quandle
$R_3^k$ and $R_3^{k-2} \times  \Z_3 [t]/ (t+1)^2$, respectively. 
$\Box$

\begin{proposition} \label{medprop}
The Galkin quandle $G(A,c)$ is medial if and only if $3A = 0$.
\end{proposition}
{\it Proof.\/}  We have seen that if $3A = 0$ then $G(A,c)$ is
Alexander and hence is medial. 
Suppose $3b \neq 0$ for some $b \in A$. Then consider the products
\begin{eqnarray*}
  X &=&  ((0,0)*(1,b))*((1,0)*(0,0) ) = (-1, b-\tau(-1)) \quad \mbox{\rm and}  \\
  Y &=&  ((0,0)*(1,0))*((1,b)*(0,0) ) = (-1,-\tau(-1) -2b) . 
  \end{eqnarray*}
Since $3b \neq 0$ we have $X \neq Y$ and so $G(A,c)$ is
not medial.
$\Box$

\begin{remark} 
{\rm 
The fact that the same condition appeared  in \ref{leftlem}, \ref{alexprop} and \ref{medprop}
is explained as follows. 
Alexander quandles are left-distributive and medial. 
It is easy to check  that for a finite Alexander quandle $(M, T)$ with $T \in {\rm Aut}(M)$,
the following are equivalent:
(1) $(M,T)$ is connected, 
(2) $(1-T)$ is an automorphism of $M$, and 
(3)  $(M,T)$ is Latin. 
It was also proved by Toyoda~\cite{Toyo} that 
a Latin quandle is Alexander if and only if it is medial. 
As  noted by Galkin, $G(\Z_5,0)$ and $G(\Z_5,1)$ are  the smallest non-medial Latin quandles
and hence the smallest non-Alexander Latin quandles.

We note that medial 
quandles are  left-distributive (by idempotency). 
We show in Theorem~\ref{leftlatinthm} (below) 
that  any left-distributive connected quandle is Latin.
This 
 implies, by Toyoda's theorem, that 
every medial connected quandle is Alexander and Latin.
The smallest Latin quandles that are not left-distributive are the Galkin
quandles of order $15$.  

It is known that the smallest left-distributive Latin quandle that is not Alexander
is of order $81$. This is due to V. D. Belousov. See, for example, \cite{Pf} or Section $5$ of
Galkin~\cite{galkin}.

}
\end{remark}

\begin{theorem} \label{leftlatinthm}
Every finite connected left-distributive connected quandle is Latin.
\end{theorem}
{\it Proof.\/}
Let $(X,*)$ be a finite, connected, and left-distributive quandle. For each $a\in X$,
 let $X_a=\{a*x:x\in X\}$.

$1^\circ$
We claim that $|X_a|=|X_b|$ for all $a,b\in X$. For any $a,y\in X$, we have
\[
|X_a|=|X_a*y|=|\{(a*x)*y:x\in X\}|=|\{(a*y)*(x*y):x\in X\}|=|X_{a*y}|.
\]
Since $X$ is connected, we have $|X_a|=|X_b|$ for all $a,b\in X$.

$2^\circ$
Fix $a\in X$. If $|X_a|=X$, by $1^\circ$, $X_b=X$ for all $b\in X$ and we are done. So assume $|X_a|<|X|$. Clearly, $(X_a,*)$ is a left-distributive quandle. Since $(X,*)$ is connected and $x\mapsto a*x$ is an onto homomorphism form $(X,*)$ to $(X_a,*)$, $(X_a,*)$ is also connected. Using induction, we may assume that $(X_a,*)$ is Latin.

$3^\circ$
For each $y\in Y$, we claim that $X_{a*y}=X_a$. In fact,
\[
\begin{split}
X_{a*y}\,&\supset (a*y)*X_a\cr
&=X_a\qquad\text{(since $X_a$ is Latin)}.
\end{split}
\]
Since $|X_{a*y}|=|X_a|$, we must have $X_{a*y}=X_a$.

$4^\circ$ Since $(X,*)$ is connected, by $3^\circ$, $X_b=X_a$ for all $b\in X$. Thus $X=\bigcup_{b\in X}X_b=X_a$, which is a contradiction.
$\Box$

\begin{proposition}\label{dualprop}
Any Galkin quandle is self-dual, that is, isomorphic to its dual.
\end{proposition}
{\it Proof.\/}
The dual quandle structure of $G(A, \tau)=G(A, c_1, c_2)$ is 
written by 
$$(x,a) \ \bar{*}\  (y, b)=(x\ \bar{*}\  y, -a + \mu(y-x)b + \tau(y-x))$$
for $(x,a),  (y,b) \in G(A, \tau)$.
Note that  $\mu(x-y)=\mu(y-x)$
and $\tau(y-x)=c_{-i}$ if $\tau(x-y)=c_i$  for any $x, y \in X$
and  $i \in \Z_3$. Hence its dual is $G(A, c_2, c_1)$.
The isomorphism is given by 
$f : \Z_3 \times A \rightarrow \Z_3 \times A$ that is defined 
 by $ f(x,a) = (-x,a)$. 
$\Box$

\begin{corollary}\label{keicor}
For any positive integer $n$, a  Galkin quandle $G (A, c_1, c_2)$ is involutory $($kei$)$
if and only if $c_1=c_2 \in A$.
\end{corollary}
{\it Proof.\/}
A quandle is a kei if and only if 
 it is the same as its dual, i.e., the identity map 
is an isomorphism between the dual quandle and itself. 
Hence this follows from Proposition~\ref{dualprop}.
$\Box$

\bigskip

A {\it good involution} \cite{K,KO} $\rho$ on a quandle $(X, *)$ is an involution
$\rho: X \rightarrow X$ 
(a map with $\rho^2=\id$)  
such that 
$x * \rho(y)=x \ \bar{*}\  y$ and $\rho(x*y)=\rho(x) * y$ for any $x, y \in X$. 
A quandle with a good involution is called a {\it symmetric} quandle.
A kei is a symmetric quandle with $\rho=\id$ (in this case $\rho$ is said to be trivial). 
Symmetric quandles have been used for unoriented knots and 
non-orientable surface-knots.

Symmetric quandles with non-trivial good involution have been hard to find.
Other than computer calculations, very few constructions have been known.
In \cite{K,KO}, non-trivial good involutions were 
defined on dihedral quandles of even order, which are not connected.
Infinitely many symmetric connected quandles were constructed in \cite{COS}
as extensions of 
odd order dihedral quandles: For each odd $2n+1$ ($n \in \Z$, $n>0$), a symmetric 
connected quandle of order $ (2n+1) 2^{2n+1}$ was given, that are not keis. 
Here we use Galkin quandles to construct more symmetric quandles.

\begin{proposition}\label{symmprop}
For any positive integer $n$,   there exists a symmetric 
connected quandle of order $6n$, that is not involutory. 
\end{proposition}
{\it Proof.\/}
We show that if  an  abelian group $A$ has an 
element $c \in A$ of order $2$, then   $G(A, c)$ is a symmetric quandle.
(A finite abelian group $A$ has an element of order $2$ if and only if $|A|$ is even.)
 Note that $G(A, c)$ is not involutory by Corollary~\ref{keicor}. 
 
Define an involution $\rho: \Z_3 \times A \rightarrow \Z_3 \times A$ 
by $\rho(x,a )=(x, a + c )$,
where  $ c \in A$ is a fixed element of order $2$ and $x \in \Z_3$, $a \in A$. 
The map $\rho$ is an involution. 
It satisfies the required conditions as we show below. 
For $x, y \in \Z_3$,  we have 
\begin{eqnarray*}
 ( x, a) * \rho( y, b)  
&=&   ( x, a)  *  (y,   b+c )  \  = \  (2y-x, -a + \mu(x-y) (b+c ) + \tau(x-y)
 , \\
 ( x,  a) \ \bar{*} \  ( y, b )  
&=& 
(  2y-x ,  - a + \mu(y-x) b  + \tau(y-x) ) ,
\end{eqnarray*} 
where the last equality follows from the proof of Proposition~\ref{dualprop}. 
If $x=y$, then $\mu(x-y)=2=\mu(y-x)$ and  $\tau(x-y)=0=\tau(y-x)$, 
and the above two terms are equal. 
If $x \neq y$,  then $\mu(x-y)=-1=\mu(y-x)$, and exactly one of 
$\tau(x-y)$ and $\tau(y-x)$ is $c$ and the other is $0$, so that the equality holds.

Next we compute 
\begin{eqnarray*}
\lefteqn{
\rho( \  ( x, a) * ( y, b )  \ )
} \\
&=& 
\rho( 2y-x , -a + \mu(x-y) b + \tau(x-y) ) \ = \ 
 (2y-x , -a +   \mu(x-y) b + \tau(x-y) + c  ), \\
 \lefteqn{
\rho  ( x, a ) * ( y, b )
}\\
 &=&
(x, a +c )* ( y, b ) 
\ = \ 
(2y-x,   -a - c + \mu (x-y) b + \tau(x-y) ) 
\end{eqnarray*}
and these are equal. 
$\Box$

\bigskip

For equations in Lemma~\ref{mutaulem}, we have the following for  $\Z_p$.

\begin{lemma}{\rm
Let $p > 3$ be a prime and let $\mu : \Z_p \rightarrow \Z$ be a function satisfying
$ 
\mu (0) = 2 $ and 
\begin{equation}\label{mu22}
 \mu (x + y) + \mu (x - y) =  \mu (x) \mu (y)
\end{equation}
for any $x, y \in \Z_p$. 
Then $\mu (x) = 2$ for all $x \in \Z_p$.
}
\end{lemma}
{\it Proof.\/}
First we note that 
we only need to 
prove 
$ \sum_{x\in\Bbb Z_p}\mu(x)\ne0.
$ 
Denote this sum by $S$. Summing Equation \eqref{mu22}
as $y$ runs over $\Bbb Z_p$, we have 
$ 
2S=S\mu(x).
$ 
So if $S\ne 0$, we have $\mu(x)=2$ for all $x\in\Bbb Z_p$.

Assume to the contrary that $S=0$.
Since $\mu(kx)\mu(x)=\mu((k+1)x)+\mu((k-1)x)$,
it is easy to see by induction that
\begin{equation}\label{11}
\mu(x)^k=\frac 12\sum_{0\le i\le k}\binom ki\mu((k-2i)x).
\end{equation}
(Here we also use the fact that $\mu(-x)=\mu(x)$, which follows from the fact that $\mu(x-y)=\mu(x)\mu(y)-\mu(x+y)$ is symmetric in $x$ and $y$.)
In particular,
\[
\mu(x)^{2p}=\frac 12\sum_{0\le i\le 2p}\binom{2p}i\mu(2(p-i)x).
\]
Since $\sum_{x\in\Bbb Z_p}\mu(x)=0$, we have
\[
\sum_{x\in\Bbb Z_p}\mu(x)^{2p}=\Bigl[2+\binom{2p}p\Bigr]p.
\]
Since $\mu(x)=\mu(\frac x2)^2-2$, we have $\mu(x)=-2,-1,2,7,\cdots$. 

{\bf Case 1.} Assume that 
there exists $0\ne x\in\Bbb Z_p$ such that $\mu(x)\ge 7$. Then
\[
\Bigl[2+\binom{2p}p\Bigr]p=\sum_{x\in\Bbb Z_p}\mu(x)^{2p}\ge 7^{2p},
\]
which is not possible.

{\bf Case 2.} Assume that $\mu(x)\in\{-2,-1,2\}$ for all $x\in\Bbb Z_p$. Let
$a_i=|\mu^{-1}(i)|$. Since $\sum_{x\in\Bbb Z_p}\mu(x)=0$ and $\sum_{x\in\Bbb Z_p}\mu(x)^3=0$, where the second equation follows from \eqref{11}, we have
\[
\begin{cases}
-2a_{-2}-a_{-1}+2a_2=0,\cr
-8a_{-2}-a_{-1}+8a_2=0.
\end{cases}
\]
So $a_{-1}=0$, i.e., $\mu(x)=\pm 2$ for all $x\in\Bbb Z_p$. Then 
\[
\sum_{x\in\Bbb Z_p}\mu(x)\equiv 2p\equiv 2\pmod 4,
\]
which is a contradiction.
$\Box$ 

\section{Knot colorings by Galkin quandles} \label{knotsec}

In this section we investigate knot colorings by Galkin quandles.
Recall from Lemma~\ref{subRnlem} 
 that any Galkin quandle $G(A, \tau)$  consists of three disjoint 
subquandles $\{x\} \times A$ for $x \in \Z_3$, 
and each is a product of dihedral quandles. 
Also any Galkin quandle has $R_3$ as a quotient. 
Thus we look at relations between colorings by dihedral quandles and those by Galkin quandles.

First we present the numbers of   $n$-colorable 
knots  with $12$ crossings  or less, 
out of $2977$ knots in the knot table from \cite{KI},
for comparison with Table~1.
These are for dihedral quandles and their products that 
may be  of interest and 
relevant for comparisons.

\vspace{-10mm}

\begin{center}
$$ 
\begin{array}{lrrlrrlrrlrlrrlrr}
R_3: & 1084 ,& &
R_5 :& 670, & &
R_7 : & 479 ,& &
R_{11} : &  285, & &
R_{15} : & 1512,  & & 
R_{17} : & 192,   \\
R_{19}:  & 159, & &
R_{21} : & 1386, & &
R_{23} : & 128, & &
R_{29}: & 97 & &
R_{31} : &  87, & &
R_{33} : &  1260. 
\end{array}
$$
\end{center}

\begin{remark}
{\rm 
We note that many rig Galkin quandles in Table~1
have the same number 
$(1084)$
of non-trivially colorable knots as the number of  
$3$-colorable knots.
We make a few observations on these Galkin quandles.

By Lemma~\ref{R3lem}, 
a Galkin quandle has $R_3$ as a subquandle if 
  $\tau=0$ or 
$3$ divides $ |A| $, and among rig Galkin quandles with the number $1084$, $17$ of them satisfy this condition. 
Hence any $3$-colorable knot is non-trivially colored by these Galkin quandles.
The converse is not necessarily  true:
$G(\Z_5, 0)$ has $\tau=0$ but has  the number $1512$.  See   Corollary~\ref{3pcor}
for more on these quandles. 

The remaining $6$  rig Galkin  quandles with the number $1084$ have 
$C[6,2]$ as a subquandle:
$$
C[12, 5], \ C[12, 9], \ C[24, 28], \ C[24, 29], \ C[24, 31], \ C[24, 33]. $$
It was conjectured \cite{COS} that 
if a knot is $3$-colorable, then it is non-trivially colored by $C[6,2]$
($\tilde{R}_3$ in their notation).
It is also seen that any non-trivial coloring by $C[6,2]$ descends to 
a non-trivial $3$-coloring via the surjection $C[6,2] \rightarrow R_3$, so if the conjecture
is true, then any knot is non-trivially colored by these quandles 
if and only if it is $3$-colorable.
See also Remarks~\ref{samenbrem} and \ref{samenontrivrem}. 
}
\end{remark}

\begin{proposition} \label{colorprop}
Let $K$ be a knot with a prime determinant  $p>3$. 
Then $K$ is non-trivially colored by a finite Galkin quandle 
$G(A, \tau)$ if and only if $p$ divides $|A|$. 
\end{proposition}
{\it Proof.\/}
Let  $K$ be a knot with the determinant that is a prime $p>3$.
By Fox's theorem~\cite{Fox}, for any prime $p$, a knot is $p$-colorable  if
and only if its determinant is divisible by $p$. 
Hence $K$ is  $p$-colorable, and not $3$-colorable.

Let $G(A, \tau)$ be any Galkin quandle and
${\cal C}: {\cal A} 
 \rightarrow G(A, \tau)$ be a coloring, where 
${\cal A}$ 
is the set of arcs of a knot diagram of $K$. 
By the surjection $r: G(A, \tau) \rightarrow R_3$, the coloring ${\cal C}$ 
induces a coloring $r  \circ {\cal C}: {\cal A} 
  \rightarrow R_3$. 
Since $K$ is not $3$-colorable, it is a trivial coloring, 
and therefore, ${\cal C} ( {\cal A} ) 
\subset r^{-1} (x)$ for some $x \in R_3$. 
The subquandle $r^{-1}(x)$ for any $x \in R_3$ is 
 an Alexander quandle 
$\{x\} \times A$ with the operation 
$(x, a)*(x, b)=(x, 2b-a)$, so that it is a product of dihedral quandles
$\{ x \} \times A=R_{q_1} \times \cdots \times  R_{q_k}$ for some positive integer $k$ and prime powers $q_j$, $j=1, \ldots, k$
(Lemma~\ref{subRnlem}).  
It is known  that the number of colorings by a product quandle 
$X_1 \times \cdots \times  X_k$
is the product of numbers of colorings by $X_i$ for $i=1, \ldots, k$. 
It is also seen that a knot is non-trivially colored by $R_{p^k}$ for a prime $p$ if and only if 
it is $p$-colorable.   
Hence 
$K$ is non-trivially colored by  $\{x\} \times A$
 if and only if one of 
$q_1, \ldots, q_k$ 
 is a power of $p$. 
$\Box$

\begin{corollary}\label{samecolcor}
For any positive integer $n$ not divisible by $3$
and any finite Galkin quandle $G(A, \tau)$, 
all $2$-bridge knots with the determinant $n$ 
have the same number of colorings by $G(A, \tau)$. 
\end{corollary}
{\it Proof.\/}
Let $K$ be a two-bridge knot with the determinant $n=p_1^{m_1} \cdots p_\ell^{m_\ell}$
 (in the prime decomposition form),
where $p_i \neq 3$ for $i=1, \ldots, \ell$,
and let $A=R_{q_1} \times \cdots \times R_{q_k}$ be the decomposition for prime powers,
as a quandle. 
By Fox's theorem~\cite{Fox}, 
for a prime $p$, $K$ is $p$-colorable if and only if $p$ divides the determinant of $K$.
Hence $K$ is $p_i$-colorable for $i=1, \ldots, \ell$, and not $3$-colorable.
By the proof of Proposition~\ref{colorprop}, the number of colorings by a Galkin quandle 
$G(A, \tau)$ of $K$ is determined by the number of colorings by the dihedral quandles 
$R_{q_j}$ 
 that are factors of $A$. 

The double branched cover $M_2(K)$ of the $3$-sphere ${\mathbb S}^3$ along 
a $2$-bridge knot $K$ is a lens space (\cite{Rolf}, for example) and its 
first homology group $H_1(M_2(K), \Z)$  is cyclic. 
If the determinant of $K$ is $n$ 
then it is 
isomorphic to $\Z_n$ 
(\cite{Lick}, for example). 
It is known \cite{Pr} that the 
number of colorings by $R_{q_j}$ 
 is equal to 
the order of 
the group
$(\Z \oplus H_1(M_2(K), \Z)) \otimes \Z_{q_j}$,  
which is determined by $n$ 
and $q_j $  
 alone. $\Box$

\begin{example} {\rm
Among knots with $8$ crossings or less,
the following sets of knots have the same numbers of colorings by all finite Galkin quandles 
from   
 Corollary~\ref{samecolcor}:
$\{ 4_1, 5_1 \}$ (determinant $5$), $\{5_2, 7_1 \}$ ($7$), 
$\{ 6_2, 7_2 \}$ ($11$), 
$\{ 6_3, 7_3, 8_1 \}$ ($13$), 
$\{ 7_5, 8_2, 8_3 \}$ ($17$), 
$\{ 7_6, 8_4 \}$ ($19$), 
$\{ 8_6, 8_7 \}$ ($23$), 
$\{ 8_8, 8_9 \}$  ($25$), 
$\{ 8_{12}, 8_{13} \}$ ($29$).
This exhausts such sets of knots up to $8$ crossings.

Computer calculations show that 
the  set of knots up to $8$ crossings with determinant $9$ is $\{6_1, 8_{20} \}$,
and these have different numbers of colorings by some Galkin quandles. 
The determinant was looked up at KnotInfo~\cite{KI}.

There  are two knots 
($7_4$ and $8_{21}$,  up to $8$ crossings) 
with determinant $15$. 
They can be distinguished by the numbers of colorings by some Galkin
quandles, according to computer calculations.
}
\end{example}

\begin{corollary} \label{3pcor}
Let $p$ be an odd prime.
Then a knot  $K$ is non-trivially colored by the Galkin quandle 
$G(\Z_p, 0)$ if and only if it is $3p$-colorable.
\end{corollary}
{\it Proof.\/} 
Suppose it is $3p$-colorable, then it is non-trivially colored by $R_{3p}$ 
which is isomorphic to 
$R_3 \times R_p$, so that it is either $3$-colorable or $p$-colorable. 
If $K$ is $3$-colorable, then 
since $G(\Z_p, 0)$ has $R_3$ as a subquandle by Lemma~\ref{R3lem}, 
$K$ is non-trivially colored by $G(\Z_p, 0)$.
If $K$ is $p$-colorable, then since   $G(\Z_p, 0)$ has $\{0 \} \times R_p$ as a subquandle by Lemma~\ref{subRnlem}, $K$ is non-trivially colored by $G(\Z_p, 0)$.

Suppose 
a knot $K$ is non-trivially colored by $G(\Z_p, 0)$, where $p$ is an odd prime.
If $K$ is $3$-colorable, then it is $3p$-colorable, and we are done.
If $K$ is not $3$-colorable, then by the proof of Proposition~\ref{colorprop}, 
$K$ is non-trivially colored by $\{x \} \times R_p$, where $x \in \Z_3$.
Hence $K$ is $p$-colorable, and so $3p$-colorable. 
$\Box$

\begin{remark} \label{samenbrem}
{\rm
According to computer calculations, 
the following sets of Galkin quandles (in the 
numbering of Table~1) 
have the same numbers of 
colorings 
for all  $2977$ knots with $12$ crossings or less. 
Thus we conjecture that it is the case for all knots.
If a Galkin quandle does not appear in the list, then it means that 
it has different numbers of colorings for some knots, comparing to other Galkin quandles.
The numbers of colorings are distinct for distinct sets listed below as well. 
\begin{eqnarray*}
& &  \{ C[6,1],  C[6,2] \},
\{ C[12,5] , C[12,6] \}, 
\{ C[12,8],  C[12,9]  \},
\{C[18,1] ,C[18,4] \}, \\ 
& & 
\{C[18,5] ,C[18,8] \}, \{C[24,27] ,C[24,28]\}, 
\{C[24,29],C[24,30] , C[24,31] \}, \\
& & 
\{ C[24,38] ,C[24,39] \},
\{ C[30,12] ,C[30,14] \},
\{C[30,13] , C[30,15] \}.
\end{eqnarray*}

} \end{remark}

\begin{remark} \label{samenontrivrem}
{\rm

In contrast to the preceding remark, if we relax the requirement of coloring the 
same number of times, and instead
consider two quandles equivalent if each colors the same knots non-trivially
(among these  $2977$  knots),
then we get the following $4$ equivalence classes. 
\begin{eqnarray*}
&& \{ C[3,1] , C[6,1] , C[6,2] , C[9,2] , C[9,6] , C[12,5] , C[12,6] , C[12,8] , C[12,9] , 
 C[18,1] , 
C[18,4] , \\ & & C[18,5] , C[18,8] ,
C[24,27] , C[24,28] , C[24,29] , C[24,30] , C[24,31] , C[24,38] , 
 C[24,39] , C[27,2] , \\ & & 
 C[27,12] , C[27,13] , C[27,23] , C[27,55] \}, \\
& & \{ C[12,7] , C[24,32] , C[24,33] \} , \\
& & \{C[15,5] , C[30,12] , C[30,14] \}, \\
 & & \{  C[15,6] , C[30,13] , C[30,15] \}. 
 \end{eqnarray*}
 
 Thus we conjecture that it is the case for all knots. 
Of these, the first family with many elements consists of quandles with 
$C[3,1]$, $C[6,1]$ or $C[6,2]$ as a subquandle. 
Hence, in fact, the conjecture about this family follows from the conjecture 
about $ \{ C[6,1],  C[6,2] \}$ in the preceding remark.

} \end{remark}

\begin{remark}{\rm
Also in contrast to 
Remark~\ref{samenbrem}, 
there exists a virtual knot $K$ (see, for example, \cite{Lou})
such that the numbers of colorings by $C[6,1]$ and $C[6,2]$ are distinct.
A  virtual knot $K$ with the following property 
was given  in \cite{COS}, Remark 4.6:
$K$ is  $3$-colorable, but does not have  a non-trivial coloring
by $C[6,2]$. 
Since $C[6,1]$ has $R_3$ as a subquandle, this virtual knot $K$ 
has a non-trivial coloring by $C[6,1]$. Hence the numbers of colorings 
by $C[6,1]$ and $C[6,2]$ are distinct for $K$. 
Thus we might conjecture that 
for any pair of non-isomorphic 
Galkin quandles, there is a virtual knot 
with different numbers of colorings. 

} \end{remark}

\begin{remark}{\rm
For any 
finite 
Galkin quandle
$G(A, \tau)$, there is a knot $K$ with a surjection $\pi_Q(K) \rightarrow G(A, \tau)$
from the fundamental quandle $\pi_Q(K)$. 
In fact, a connected sum of trefoil can be taken as $K$ as follows
(see, for example, \cite{Rolf} for connected sum). 

First we take a set of generators of $G(A, \tau)$ as follows.
Let $A=\Z_{n_1}\times \cdots \times \Z_{n_k}$, where 
$k, n_1, \ldots , n_k$ are positive integers such that $n_i$ divides $n_{i+1} $
for $i=1, \ldots, k$. 
Let $S=\{  (x, e_i) \ | \ x \in \Z_3, i=0, \ldots, k \}$, 
where $e_0=0 \in A$ and $e_i \in A$ ($i=1, \ldots, k$) is an elementary vector
$[0, \ldots, 0, 1, 0, \ldots, 0] \in \Z_{n_1} \times \cdots \times \Z_{n_k}$
with a single $1$ at 
the 
$i$-th position.
Note that $R_n$ is generated by $0, 1$ as 
$0*1=2$, $1*2=3$, and inductively, $i*(i+1)=i+2$ for $i=0, \ldots, n-2$. 
Since $\{x \} \times A$ is isomorphic to a product of dihedral quandles
for each $x \in \Z_3$, $S$ generates $G(A, \tau)$.

For a $2$-string braid $\sigma_1^3$  whose closure is trefoil (see Figure~\ref{tref}), 
we note that if $x \neq y \in \Z_3$, then for any 
$a, b \in A$, the pair of colors $(x,a), (y, b) \in G(A, \tau)$ at top arcs extends to the bottom,
i.e., the bottom arcs receive the same pair. This can be computed directly.

For copies of trefoil, we assign
 pairs $[(0, e_0), (x, e_i)]$ as colors 
where $x=1,2$ and $i=0, \ldots k$, 
and take connected sums on the portion of the arcs with the common color 
$(0, e_0)$. 
Further we take pairs $[(0, e_j), (1, e_0)]$ for $j=1, \ldots, k$, for example, 
and take connected sum on the arcs with the common color $(1, e_0)$,
to obtain a connected sum of trefoil with all elements of $S$ used as colors. 
Such a coloring gives rise  to a quandle homomorphism 
$\pi_Q(K)  \rightarrow G(A, \tau)$ whose image contains generators $S$,
hence defines a surjective homomorphism. 
}
\end{remark}

\end{document}